\documentclass[a4paper,11pt]{article}
\usepackage{latexsym,amsfonts,amsmath}
\usepackage{amssymb}
\usepackage{color}
\usepackage[mathscr]{eucal}
\usepackage{enumerate}

\def\1{\mathbf{1}}

\def\NN{\mathbb{N}}

\def\RR{\mathbb{R}}

\newcommand{\widebar}[1]{\overline{#1}}

\def\union{\mathop{\cup}}

\def\liminf{\mathop{\underline{\lim}}}
\def\limsup{\mathop{\overline{\lim}}}
\def\ds{\displaystyle}

\def\nl{\mbox{} \newline }

\newcounter{hypot}

    \newenvironment{hypot}{\begin{list}
      {\hspace{\labelsep}\bfseries Assumption \Alph{hypot}.}
      {\leftmargin=0pt
       \labelwidth=0cm
       \refstepcounter{hypot}
       \def\makelabel##1{##1}}}{\end{list}}

\newcounter{assump}

\newenvironment{assump}{\begin{list}
      {\hspace{\labelsep} (\Alph{hypot}.\arabic{assump})}
      {\leftmargin=0pt
       \labelwidth=0cm
       \usecounter{assump}
       }}{\end{list}}

\begin{document}
\newtheorem{theorem}{Theorem}[section]
\newtheorem{proposition}[theorem]{Proposition}
\newtheorem{lemma}[theorem]{Lemma}
\newtheorem{corollary}[theorem]{Corollary}
\newtheorem{definition}[theorem]{Definition}
\newtheorem{remark}[theorem]{Remark}
\newtheorem{conjecture}[theorem]{Conjecture}
\newtheorem{assumption}[theorem]{Assumption}

\bibliographystyle{plain}

\title{On  the expected total reward with unbounded returns for Markov decision processes}

\author{ \mbox{ }
F. Dufour \\
\small Institut Polytechnique de Bordeaux \\
\small INRIA Bordeaux Sud Ouest, Team: CQFD\\
\small IMB, Institut de Math\'ematiques de Bordeaux, Universit\'e de Bordeaux, France\\
\small e-mail: francois.dufour@math.u-bordeaux.fr
\and
A. Genadot \\
\small IMB, Institut de Math\'ematiques de Bordeaux, Universit\'e de Bordeaux, France\\
\small INRIA Bordeaux Sud Ouest, Team: CQFD\\
\small e-mail: alexandre.genadot@math.u-bordeaux.fr
}

\date{}

\maketitle

\begin{abstract} We consider a discrete-time Markov decision process with Borel state and action spaces. The performance criterion is to maximize a total expected {utility determined by unbounded return function.} 
It is shown the existence of optimal strategies under general conditions allowing the reward function to be unbounded both from above and below and the action sets available at each step to the decision maker to be not  necessarily compact.
To deal with unbounded reward functions, a new characterization for the weak convergence of probability measures is derived. 
Our results are illustrated by examples.
\end{abstract}

{\small
}

\section{Introduction}
\label{sec-1}
In this paper, our objective is to provide sufficient conditions for the existence of optimal strategies in dynamic programming decision models under the expected total reward criterion. The model under consideration is rather general since the reward function may be unbounded both from above and below and the action sets available at each step to the decision maker may not be necessarily compact.

Typically, here is an example of model we are able to handle. Roughly speaking, the state space is given by $\mathbf{X}=[0,1]$, the action space is {$\mathbf{A}=\{1,2,\ldots\}$} and the reward function $r$ satisfies
$$r(x,a)=\begin{cases}
0 & \text{ if $(x,a) \in \{0\}\times {\mathbf{A}}$} \\
\frac{1}{a^{2}} \frac{1}{\sqrt{x}}I_{]0,1/2]}(x)-a \frac{1}{\sqrt{1-x}}I_{]1/2,1[}(x) & \text{ if  $ (x,a) \in ]0,1[\times {\mathbf{A}} $} \\
-\infty, & \text{ if $(x,a) \in \{1\} \times {\mathbf{A}}$.}
\end{cases}$$
Observe that this reward function is not upper semicontinuous nor bounded from above. Moreover, it can actually takes the $-\infty$ value. This example will be described in details in Section \ref{sec-examples}. As far as we know, such example does not satisfy the standard conditions of the literature see, for instance, the references \cite{balder89,balder92,nowak88,schal75,schal79}.

{Several approaches have been proposed in the literature to study the existence of optimal strategies for discrete-time Markov decision processes under the total expected utility criterion with unbounded reward function. A possible method is based on the analysis of the so-called dynamic programming equation, see for example
\cite{bertsekas78,hinderer70,schal75} for general results in this direction.
In addition, the total reward criterion can be seen as a special case of the so-called expected utility criteria as studied for example in \cite{kertz79}.}
In our paper, we use a different approach following the line of the works developed in \cite{balder89,balder92,nowak88,schal75,schal79} that, roughly speaking, consists in finding a suitable topology on the set of strategic probability measures {to ensure that this set is compact and the expected reward functional is semicontinuous.} This technique is quite classical in the literature. In the context of compact action sets, it has been studied by many authors, see \textit{e.g.} \cite{balder89,nowak88,schal75,schal79}. In \cite{balder92}, this method has been generalized to the case of possible non-compact action sets by introducing the so-called strong coercivity condition on the reward function.
A key condition used in all the aforementioned works \cite{balder89,balder92,nowak88,schal75,schal79} is to consider the reward function being  bounded from above (or equivalently, bounded from below for a cost function).
However, it has been emphasized in \cite{jaskiewicz14,jaskiewicz11,matkowski11} (and the references therein) that many real models do not satisfy such property, in particular in economy where the utility function may be logarithmic and may take the value $-\infty$ for some states (see example 2 in \cite{jaskiewicz11}). 
Recently, new sets of conditions have been studied in \cite{jaskiewicz14,jaskiewicz11,matkowski11} to deal with Markov decision processes {and in \cite{jaskiewicz11b}
for min-max games} generalizing the so-called weighted norm approach due to Wessels \cite{wessels77}. Roughly speaking the authors provide existence results based on the analysis of the Bellman optimality equation.
In this paper, \textit{red}{we used another approach} and follow the line developed in \cite{balder92} to show the existence of optimal strategies in the general framework of unbounded reward function and so generalizing the results of \cite{balder92}.
Moreover, it will be shown that the assumptions proposed in \cite{jaskiewicz11} satisfy our new set of hypotheses.

In order to deal with unbounded reward functions, we have derived a new characterization for the weak convergence of probability measures. We show that the determining class of test functions for weak convergence, usually the set of continuous and bounded functionals, can be relaxed to the set of not everywhere bounded and upper semicontinuous functions. Such a relaxation on the boundedness condition have also been considered in \cite{zapala08} but in keeping a continuity hypotheses. 

\bigskip

The rest of the paper is organized as follows. In Section \ref{sec-model}, we define the control model.
Section \ref{sec-portmanteau} studies weak convergence of probability measures with not everywhere bounded and semicontinuous test functions.
In Section \ref{sec-results} we state our assumptions and establish the existence of optimal policies.
Finally, Section \ref{sec-examples} is dedicated to the presentation of examples illustrating our results.
\section{Description of the control problem}
\label{sec-model}
\setcounter{equation}{0}

\subsection{Basic notations}

First of all, we introduce the following notations and terminology. 

We write $\NN$ for the set of integers, that is, $\NN=\{0,1,2,\ldots\}$,  $\NN^*=\NN\setminus\{0\}$ for the set of positive integers, $\RR$ for the set of real numbers, $\widebar{\RR}=\RR\cup\{-\infty,+\infty\}$ for the set of extended real numbers.

Given a topological space $Z$, we say that a function $f:Z\rightarrow\widebar{\RR}$ is upper semicontinuous (respectively, upper semicompact) on $Z$ if $\{z\in Z: f(z)\geq \beta\}$ is a closed (respectively, relatively compact) set of $Z$ for every
$\beta\in\RR$. 
By $\widebar{\boldsymbol{\mathcal{C}}}(Z)$ we will denote the family of real-valued functions on $Z$ which are
bounded and continuous.

Recall from \cite{balder92} that, for two metric spaces $Y$ and $Z$ and a subset $B$ of $Y\times Z$, a function $u : B \to [-\infty,+\infty)$ is said to be strongly coercive on $B$ if for every sequence $\{(y_k,z_k)\}_{k\in \NN}$ in $B$, such that $y_k\to y^{*}$ for some $y^{*}\in Y$ and $\beta =\limsup_{k} u(y_k,z_k)>-\infty$, there exists a subsequence $\{z_{k_j}\}$ of $\{z_k\}$ and $z^*\in Z$ such that $z_{k_j}\to z^*$, $(y^{*},z^{*})\in B$ and $u(y^{*},z^*)\geq\beta$. As a direct consequence of the definition, strongly coercive functions are upper semicontinuous. However, it is not hard to exhibit upper semicontinuous fonctions which are not strongly coercives.
This definition of strongly coercivity is consistent with the notion of coercivity in the literature. Indeed, if $u : Y\times Z\to [-\infty,+\infty)$ is strongly coercive with $(Z,\|\cdot\|)$ a normed vector space, then, for any $y\in Y$, we have
$\ds \lim_{\|z\|\to\infty} u(y,z)=-\infty$.

We use the symbol $f^{+}$ (respectively $f^{-}$) to denote the positive part (respectively, negative part) of a function
$f:Z\rightarrow\widebar{\RR}$.

The Borel $\sigma$-algebra of $Z$ is denoted by $\mathfrak{B}(Z)$. The set of probability measures on $(Z,\mathfrak{B}(Z))$ is denoted by $\boldsymbol{\mathcal{P}}(Z)$. In this work, $\boldsymbol{\mathcal{P}}(Z)$ will be considered as a topological space equipped with the weak topology. If $\mu$ is a measure on $Z_1\times Z_2$ then $\mu_{| Z_{1}}$ denotes the marginal of the measure $\mu$ on $Z_1$.

A Borel subset of a complete and separable metric space is called Borel space. If $Z_1$ and $Z_2$ are two Borel spaces and $Q$ is a stochastic kernel on $Z_2$ given $Z_1$, then, for a function $v:Z_2\rightarrow\widebar{\RR}$, we define $Qv:Z_1\rightarrow\widebar{\RR}$ as
$$Qv(z_1):=\int_{Z_2}v^{+}(z_2)Q(dz_2|z_1)-\int_{Z_2}v^{-}(z_2)Q(dz_2|z_1),$$
provided that {one of the two terms in the right member of the previous equation is finite.}
For a measure $\mu$ on $Z_1$, we denote by $\mu Q$ the measure $\ds \int_{Z_1} Q(\cdot|z_1)\mu(dz_1)$ on $Z_2$. 

\subsection{The control model.}

Let us consider the standard discrete-time non-stationary control model:
\begin{eqnarray*}
\big(\{\mathbf{X}_{t}\}_{t\in \NN^{*}},\{\mathbf{A}_{t}\}_{t\in \NN^{*}},\{\Psi_{t}\}_{t\in \NN^{*}},\{Q_{t}\}_{t\in \NN^{*}},\{r_{t}\}_{t\in \NN^{*}},\nu\big) \label{5tuple}
\end{eqnarray*}
consisting of:
\begin{enumerate}
\item [(a)] A sequence of Borel spaces $\{\mathbf{X}_{t}\}_{t\in \NN^{*}}$ where $\mathbf{X}_{t}$ is the state space at time $t\in\NN^{*}$.
\item [(b)] A sequence of Borel spaces $\{\mathbf{A}_{t}\}_{t\in \NN^{*}}$ where $\mathbf{A}_{t}$ represents the control or action set at time $t$.
For notational convenience, we introduce recursively the set $\mathbf{H}_{t}$ of histories up to time $t$ by defining $\mathbf{H}_{1}=\mathbf{X}_{1}$ and $\mathbf{H}_{t+1}=\mathbf{H}_{t}\times \mathbf{A}_{t}\times \mathbf{X}_{t+1}$ for $t\geq 1$.
The set of histories will be denoted by $\mathbf{H}_{\infty}=\prod_{t\in \NN^{*}} \mathbf{X}_{t}\times \mathbf{A}_{t}$.
\item [(c)] A sequence of multifunctions $\{\Psi_{t}\}_{t\in \NN^{*}}$ defined recursively by setting $\Psi_{1} \colon  \mathbf{X}_{1} \to 2^{\mathbf{A}_{1}}  \backslash \{\emptyset\}$ and 
$\Psi_{t} : \boldsymbol{\mathcal{K}}_{t-1} \times \mathbf{X}_{t} \rightarrow 2^{\mathbf{A}_{t}} \backslash \{\emptyset\}$ for $t\geq 2$ where $\boldsymbol{\mathcal{K}}_{t}\subset\mathbf{H}_{t}\times \mathbf{A_{t}}$ denotes the graph of $\Psi_{t}$. It is assumed that $\boldsymbol{\mathcal{K}}_{t}\in\mathfrak{B}(\mathbf{H}_{t}\times \mathbf{A}_{t})$.
The set $\Psi_{t}(h_{t})$ for any $h_{t}\in \boldsymbol{\mathcal{K}}_{t-1} \times \mathbf{X}_{t}$ (by a slight abuse of notation $\boldsymbol{\mathcal{K}}_{0}\times \mathbf{X}_{1}$ means $\mathbf{X}_{1}$) is the set of available actions the decision maker can choose knowing the admissible history $h_{t}$ up to time $t$. By $\boldsymbol{\mathcal{H}}_{t}$ we denote the set of admissible histories up to time $t$, that is $\boldsymbol{\mathcal{H}}_{t}=\boldsymbol{\mathcal{K}}_{t-1} \times \mathbf{X}_{t}$ for $t\in \NN^{*}$.
\item [(d)] A sequence of stochastic kernels $\{Q_{t}\}_{t\in \NN^{*}}$ on $\mathbf{X}_{t+1}$ given $\boldsymbol{\mathcal{K}}_{t}$, which stands for the transition probability function from time $t$ to time $t+1$.
\item [(d)] The reward function $r_{t}:\boldsymbol{\mathcal{K}}_{t}\rightarrow [-\infty,+\infty[$ at step $t$ for $t\in \NN^{*}$.
\item [(f)] Finally, a probability measure $\nu$ on $(\mathbf{X}_{1},\mathfrak{B}(\mathbf{X}_{1}))$ describing the initial distribution of the process. 
\end{enumerate}
A control policy (a policy, for short) is a sequence $\pi=\{\pi_{t}\}_{t\in\NN}$ of stochastic kernels $\pi_{t}$ on~$\mathbf{A}_{t}$ given $\boldsymbol{\mathcal{H}}_{t}$ such that $\pi_{t}(\Psi_{t}(h_{t})|h_{t})=1$ for any $h_{t}\in \boldsymbol{\mathcal{H}}_{t}$.
Let $\Pi$ be the set of all policies.

To state the optimal control problem we are concerned with, we introduce the canonical space $(\Omega,\mathcal{F})$ consisting of the set of sample paths $\Omega=\mathbf{H}_{\infty}$
and the associated product $\sigma$-algebra $\mathcal{F}$. The projection from $\Omega$ to the $t$-th state space and the $t$-th action space are denoted by $X_{t}$ and $A_{t}$. That is, for
\[\omega=(y_{1},b_{1},\ldots,y_{t},b_{t}\ldots)\in \Omega \quad\hbox{we have}\quad X_{t}(\omega)=y_{t}\;\;\hbox{and}\;\; A_{t}(\omega)=b_{t}\] for $t\in \NN^{*}$.
Consequently, $\{X_{t}\}_{t\in \NN^{*}}$ is the state process and $\{A_{t}\}_{t\in \NN^{*}}$ is the control process.
For notational convenience, we will write $\mathcal{H}_{t}=\sigma\{X_{1},A_{1},\ldots,X_{t-1},A_{t-1},X_{t}\}$ for $t\geq 2$ and $\mathcal{H}_{1}=\sigma\{X_{1}\}$.
It is a well known result that for every policy $\pi \in \Pi$ and any initial probability measure $\nu$ on $(\mathbf{X}_1,\mathfrak{B}(\mathbf{X}_1))$ there exists a unique probability measure $\mathbb{P}^{\pi}$ on $(\Omega,\mathcal{F})$
such that the marginal of $\mathbb{P}^{\pi}$ on $\boldsymbol{\mathcal{H}}_{t}\times\mathbf{A}_{t}$ denoted by $\mathbb{P}^{\pi}_{| \boldsymbol{\mathcal{H}}_{t}\times \mathbf{A}_{t}}$ satisfies 
$\mathbb{P}^{\pi}_{| \boldsymbol{\mathcal{H}}_{t}\times \mathbf{A}_{t}}(\boldsymbol{\mathcal{K}}_{t})=1$, for any $t\in \NN^{*}$ and
$$\mathbb{P}^{\pi}(X_{1}\in B)=\nu(B), \quad \text{ for } B\in \mathfrak{B}(\mathbf{X}_{1}),$$
$$\mathbb{P}^{\pi}(X_{t+1}\in C|\mathcal{H}_{t}\vee\sigma\{A_{t}\})=Q_{t}(C|X_{1},A_{1},\ldots,X_{t-1},A_{t-1},X_{t},A_{t}) \quad \text{ for } C\in \mathfrak{B}(\mathbf{X}_{t+1}),$$ 
$$ \mathbb{P}^{\pi}(A_{t}\in D| \mathcal{H}_{t})=\pi_{t}(D|X_{1},A_{1},\ldots,X_{t-1},A_{t-1},X_{t}) \quad \text{ for } D\in \mathfrak{B}(\mathbf{A}_{t}),$$
$\mathbb{P}^{\pi}-a.s.$, for any $t\in\NN^{*}$.

We refer to $\mathbb{P}^{\pi}$ as the \textit{strategic probability measure} generated by the policy $\pi$. Observe that we have chosen to drop in the notation of $\mathbb{P}^{\pi}$ its dependence  with respect to the initial distribution $\nu\in \boldsymbol{\mathcal{P}}(\mathbf{X}_1)$. The expectation with respect to  $\mathbb{P}^{\pi}$ is denoted by $\mathbb{E}^{\pi}$.
Let $\boldsymbol{\mathcal{S}}$ be the family of strategic probability measures, that is,
\begin{equation}\label{eq-def-P-nu}
\boldsymbol{\mathcal{S}}=\{\mathbb{P}^\pi:\pi\in\Pi\}\subseteq\boldsymbol{\mathcal{P}}(\Omega).
\end{equation}

\paragraph{Statement of the control problem.}
\nl
The expected reward functional $\mathcal{J} \colon \boldsymbol{\mathcal{S}} \to [-\infty,+\infty[$ is defined by
{\begin{eqnarray*}
\mathcal{J}(\mathbb{P}^\pi)=\sum_{t=1}^{\infty} \int_{\boldsymbol{\mathcal{K}}_{t}} r_{t}^{+}(H_{t},A_{t}) d\mathbb{P}^\pi_{|\mathbf{H}_{t}\times \mathbf{A}_{t}} - \sum_{t=1}^{\infty} \int_{\boldsymbol{\mathcal{K}}_{t}} r_{t}^{-}(H_{t},A_{t}) d\mathbb{P}^\pi_{|\mathbf{H}_{t}\times \mathbf{A}_{t}}
\end{eqnarray*}
where $\mathbb{P}^\pi\in \boldsymbol{\mathcal{S}}$} and by convention $(+\infty)-(+\infty)=-\infty$.

The optimal control problem we consider consists in maximizing the expected reward $\mathcal{J}$ over the set of strategic probability measures $\boldsymbol{\mathcal{S}}$.

\section{Weak convergence with not everywhere bounded and semicontinuous functions}
\label{sec-portmanteau}

By definition, a sequence of probability measures $(\mu_n)$ on a metric space $Y$ converges weakly towards a measure $\mu$ if for any $u\in \bar{\boldsymbol{\mathcal{C}}}(X)$,
$$
\lim_{n\to\infty} \int_Y u d\mu_n=\int_Yud\mu.
$$
It is well known that by Alexandroff's theorem, also known as Portmanteau's theorem, see \cite[Corollary 8.2.5]{bogachev07}, the condition of continuity for the class of {test functions $u$} can be relaxed to semicontinuity. Weak convergence of $\mu_n$ towards $\mu$ is thus equivalent to the fact that
$$
\limsup_{n\to\infty} \int_Y u d\mu_n\leq\int_Yud\mu
$$
for any upper semicontinuous {function $u$} which is bounded from above.
The boundedness condition can itself be relaxed to some kind of uniform integrability for the test functions, see \textit{e.g.} \cite{zapala08}. For instance, let us rewrite \cite[Theorem 2]{zapala08} in our setting.

\begin{theorem}[\cite{zapala08}]
\label{thm:zapala}
A sequence $(\mu_{n})$ of probability measures on a metric space $Y$ converges weakly towards a probability measure $\mu$ if and only if, for any continuous {function $u$} which is asymptotically uniformly integrable, that is
$$
\forall\epsilon>0,~\exists n_\epsilon\in\NN,~C_\epsilon>0,~\forall n\geq n_\epsilon,\quad\int_{\{|u|\geq C_\epsilon\}} |u| d\mu_n <\epsilon,
$$
we have
$$
\lim_{n\to\infty} \int_Y u d\mu_n=\int_Y ud\mu.
$$
\end{theorem}

We propose here an approach allowing to relax both everywhere semicontinuity and boundedness conditions.

\begin{theorem}
\label{thm:newhanger}
A sequence $(\mu_{n})$ of probability measures on a metric space $Y$ converges weakly towards a probability measure $\mu$ if and only if, for any function $u:Y\rightarrow [-\infty,+\infty[$ satisfying the following conditions:
\begin{itemize}
\item $u^+$ is integrable with respect to $\mu$:
\begin{align}
\label{Hyp-lem-tight-lsc2}
\int_{Y} u^{+} d\mu <+\infty;
\end{align}
\item for any $\epsilon>0$,  there exists a closed subset $Y_{\epsilon}$ of $Y$ satisfying
\begin{align}
\label{Hyp-lem-tight-lsc1}
\sup_{n} \int_{Y\setminus Y_{\epsilon}} [u^{+}\vee 1] d\mu_{n} <\epsilon;
\end{align}
\item the restriction of $u$ on $Y_{\epsilon}$ is upper semicontinuous and bounded above;
\end{itemize}
we have, 
\begin{align}
\limsup_{n\rightarrow \infty} \int_{Y} u d\mu_{n} \leq \int_{Y} u d\mu.
\end{align}
\end{theorem}
\textbf{Proof:} The \textit{only if} part is obvious from Alexandroff's Theorem. Let us consider the \textit{if} part. From (\ref{Hyp-lem-tight-lsc1}), we have  $\ds \sup_{n}\int_{Y} u^{+} d\mu_{n} < + \infty $
and so, $\ds \int_{Y} u d\mu_{n}$ is well defined for any $n\in \NN$.
For $\epsilon>0$, consider $Y_{\epsilon}$ satisfying the hypotheses.
Write $u_{\epsilon}(y) = u(y)\wedge M_{\epsilon}$ where $M_{\epsilon}=\sup_{Y_{\epsilon}} u^{+}$.
We have
\begin{align*}
\int_{Y} u d\mu_{n} = \int_{\{u^{+}\leq M_{\epsilon} \}} u_{\epsilon} d\mu_{n} + \int_{\{u^{+} > M_{\epsilon} \}} u^{+} d\mu_{n}.
\end{align*}
Observe now that $\{u^{+} > M_{\epsilon} \} \subset Y\setminus Y_{\epsilon}$ and $u_{\epsilon}\equiv M_{\epsilon}$ on $\{u^{+} > M_{\epsilon} \}$
showing 
\begin{align*}
\int_{Y} u d\mu_{n}  \leq \int_{Y} u_{\epsilon} d\mu_{n} + \int_{Y\setminus Y_{\epsilon}} u^{+} d\mu_{n}.
\end{align*}
Therefore, equation (\ref{Hyp-lem-tight-lsc1}) implies 
\begin{align*}
\int_{Y} u d\mu_{n}  \leq \int_{Y} u_{\epsilon} d\mu_{n} + \epsilon
\end{align*}
for any $\epsilon>0$. Clearly, by using the hypotheses, {for any $\eta>0$} there exists a closed subset $Z_{\eta}$ of $Y$ such that
the restriction of $u_{\epsilon}$ on $Z_{\eta}$ is upper semicontinuous and $\sup_{n} \mu_{n}(Y\setminus Z_{\eta}) <\eta$.
Moreover, $u_{\epsilon}$ is clearly bounded above on $Y$ and so $u_{\epsilon}$ satisfies the conditions of Lemma 2.4 in \cite{balder89} whose proof is detailed in Appendix \ref{app:balder}, implying
$\ds \limsup_{n\rightarrow \infty} \int_{Y} u_{\epsilon} d\mu_{n}\leq \int_{Y} u_{\epsilon} d\mu$. Therefore,
$$\ds \limsup_{n\rightarrow \infty} \int_{Y} u d\mu_{n}  \leq \int_{Y} u_{\epsilon} d\mu + \epsilon.$$
Now, it follows from (\ref{Hyp-lem-tight-lsc2}) that $\ds \int_{Y} u d\mu$ is well defined and since $u_{\epsilon}\leq u$ we get
\begin{align*}
\limsup_{n\rightarrow \infty} \int_{Y} u d\mu_{n}  \leq \int_{Y} u d\mu + \epsilon
\end{align*}
for any $\epsilon>0$, showing the result.
\hfill$\Box$
\bigskip

Clearly, the condition (\ref{Hyp-lem-tight-lsc1}) may be relaxed to the existence, for any $\epsilon>0$, of a closed subset $Y_\epsilon$ and an integer $n_\epsilon$ such that
$$
\sup_{n\geq n_\epsilon} \int_{Y\setminus Y_{\epsilon}} [u^{+}\vee 1] d\mu_{n} <\epsilon.
$$
Then, it is not hard to see that the latter condition plus the fact that $u$ is upper semicontinuous and bounded from above on $Y_\epsilon$ implies the asymptotic uniform integrability of $u$, and the latter condition is thus somehow stronger in this sense than asymptotic uniform integrability. But in our case, $u$ is semicontinuous only on a subset $Y_\epsilon$, making, in this sense, the set of conditions of Theorem \ref{thm:newhanger} weaker than the conditions in Theorem \ref{thm:zapala}.

\bigskip

\noindent
The following result is a direct consequence of Theorem \ref{thm:newhanger}.
\begin{corollary}
\label{lem-tight-usc}
Let $Y$ be a metric space. Consider a subset $\cal P$ of $\boldsymbol{\mathcal{P}}(Y)$ and a function $u$ on $Y$ satisfying
\begin{itemize}
\item for any $\epsilon>0$, there exists a closed subset $Y_\epsilon$ of $Y$ such that
$$
\sup_{\mathbb{P}\in\cal P}\int_{Y\setminus Y_\epsilon} (u^{+}\vee 1) d \mathbb{P} <\epsilon,
$$
\item the restriction of $u$ to $Y_\epsilon$ is upper semicontinuous and bounded above.
\end{itemize}
Then, the function $\ds \mathbb{P} \to \int_Y u d\mathbb{P}$ defined on $\cal P$ is upper semicontinuous and bounded above.
\end{corollary}

\section{Existence result under general conditions}
\label{sec-results}

We start this section with the introduction and a discussion of the assumptions under consideration in this work.
Then, we will prove the existence of an optimal control strategy for the model presented in section \ref{sec-model}.

\begin{hypot}
\item \mbox{ }
\label{Regular-Hypotheses}
\begin{assump}
\item \label{Condition-C} The following condition holds:
\[ \lim_{n\rightarrow\infty} \sup_{p\geq n} \sup_{\pi\in\Pi}  \bigg[ \sum_{t=n}^{p} \mathbb{E}^\pi \big[ r_t(H_{t},A_{t})\big] \bigg]^{+}= 0.\]
\item \label{Continuity-transition-kernel} For any $t\in \NN^{*}$, $g\in \widebar{\boldsymbol{\mathcal{C}}}(\mathbf{X}_{t+1})$ and $\epsilon>0$ there exists $C_{\epsilon}$
a closed subset of $\boldsymbol{\mathcal{K}}_{t}$ satisfying
\begin{eqnarray*}
\sup_{\pi\in\Pi} \mathbb{E}^\pi\big[ \mathbf{I}_{\boldsymbol{\mathcal{K}}_{t}\setminus C_{\epsilon}}(H_{t},A_{t}) \big]<\epsilon
\end{eqnarray*}
and such that the real-valued mapping defined on $C_{\epsilon}$ by
$\ds (x_{1},a_{1},\ldots,x_{t},a_{t})\rightarrow Q_tg(x_{1},a_{1},\ldots,x_{t},a_{t})$ is continuous.
\end{assump}
\end{hypot}

\begin{hypot}
\label{Reward+Multifunction}
\item For any $t\in \NN^{*}$, $\epsilon>0$ there exists $K_{\epsilon}$, a closed subset of $\boldsymbol{\mathcal{H}}_{t}$, satisfying
\begin{eqnarray}
\label{w-tight-condition}
\sup_{\pi\in\Pi} \mathbb{E}^\pi\Big[ \mathbf{I}_{\boldsymbol{\mathcal{H}}_{t}\setminus K_{\epsilon}}(H_{t}) \big[1\vee r_t^{+}(H_{t},A_{t}) \big] \Big]<\epsilon
\end{eqnarray}
and such that the restriction of $r_{t}$ to $[K_{\epsilon} \times \mathbf{A}_{t}]\cap\boldsymbol{\mathcal{K}}_{t}$ is strongly coercive and bounded above.
\end{hypot}

\begin{remark}
\label{Hyp-discussion}
Discussion of the hypotheses.
\begin{enumerate}[(a)]
\item\label{Schal-condition}
Assumption \ref{Condition-C} is the so-called Condition (C) in Sch\"al's papers \cite{schal75b,schal75}. It is slightly weaker than Condition (A3) in Balder's result \cite{balder92}, {see also the discussion in \cite[Remark 6.11]{kertz79}}.
\item\label{Balder-condition} {Assumption \ref{Continuity-transition-kernel} is a standard condition, see for example Condition (C2) in \cite{balder89} and also Condition (A1) in \cite{balder92}. 
Specific conditions that can be expressed in terms of the primitive data of the model and implying Assumption \ref{Continuity-transition-kernel} are presented in \cite[Section 3]{balder89}.}
\item\label{New-condition-1} The condition (\ref{w-tight-condition}) in Assumption \ref{Reward+Multifunction} is new and generalizes condition (A2) in \cite{balder92} to the case where the reward function may not be bounded above. 
Observe that in \cite{balder89,balder92,schal75b,schal75} the reward functions are bounded above.
In \cite{jaskiewicz11}, the authors studied a discounted Markov decision process on general state and action spaces with possibly unbounded reward function with application to economic models.
We will show in Section \ref{example-Nowak} that the approach presented in \cite{jaskiewicz11} can be easily embedded into our framework.
Moreover, our condition (\ref{w-tight-condition}) incorporates the case of history-dependent action spaces contrary to the framework discussed in \cite{balder92}.
\item\label{New-condition-1-sufficiency}
{It can be shown easily that the following set of conditions implies that Assumption \ref{Reward+Multifunction} holds. There are written explicitly in terms of the parameters of the model 
which makes them easier to check than Assumption \ref{Reward+Multifunction}, as outlined in Section \ref{example-general} through an example:
\begin{enumerate}[(i)]
\item For any $t\in \NN^{*}$, there exists a $\RR_{+}\union\{+\infty\}$-valued measurable mapping $\Phi_{t}$ defined on $\boldsymbol{\mathcal{H}}_{t}$ satisfying
\begin{eqnarray}
 \sup_{a_{t}\in \Psi_{t}(h_{t})} \big[1\vee r_{t}^{+}(h_{t},a_{t}) \big] \leq \Phi_{t}(h_{t}).
\end{eqnarray}
\item For any $\epsilon>0$, there exists $K_{\epsilon}$ a closed subset of $\mathbf{X}_{1}$, satisfying
\begin{align}
\label{sufficient-tight-1-condition}
\int_{\mathbf{X}_{1}\setminus K_{\epsilon}}  \Phi_{1}(x_{1})  \nu(dx_{1}) < \epsilon
\end{align}
and such that the restriction of $r_{1}$ to $[K_{\epsilon} \times \mathbf{A}_{1}]\cap\boldsymbol{\mathcal{K}}_{1}$
is strongly coercive and bounded above.
\item For any $t\in \NN^{*}$ and $\epsilon>0$, there exists $K_{\epsilon}$ a closed subset of $\boldsymbol{\mathcal{H}}_{t+1}$, satisfying
\begin{align}
\label{sufficient-tight-2-condition}
\sup_{(h_{t},a_{t})\in \boldsymbol{\mathcal{K}}_{t}} \int_{\mathbf{X}_{t+1}}  \mathbf{I}_{\boldsymbol{\mathcal{H}}_{t+1}\setminus K_{\epsilon}} (h_{t},a_{t},x_{t+1})
\Phi_{t+1}(h_{t},a_{t},x_{t+1})
Q_{t}(dx_{t+1}|h_{t},a_{t}) <\epsilon
\end{align}
and such that the restriction of $r_{t+1}$ to $[K_{\epsilon} \times \mathbf{A}_{t+1}]\cap\boldsymbol{\mathcal{K}}_{t+1}$
is strongly coercive and bounded above.
\end{enumerate}}
\item {Observe that Assumption \ref{Continuity-transition-kernel} is related to the family of transition kernels $\{Q_{t}\}_{t\in \NN^{*}}$ and states roughly speaking that $Q_{t}$ is weakly continuous on a closed subset of
$\boldsymbol{\mathcal{K}}_{t}$ for $t\in \NN^{*}$ while Assumption \ref{Reward+Multifunction} is associated with the reward functions $\{r_{t}\}_{t\in \NN^{*}}$ and imposes that 
$r_{t}$ is strongly coercive on a closed subset of $\boldsymbol{\mathcal{H}}_{t}$ for $t\in \NN^{*}$.
We would like to emphasize that the closed sets involved in Assumption \ref{Continuity-transition-kernel} and \ref{Reward+Multifunction} are different by definition since in Assumption \ref{Continuity-transition-kernel} it is a closed subset of
$\boldsymbol{\mathcal{K}}_{t}$ while in Assumption \ref{Reward+Multifunction} it is a closed subset of $\boldsymbol{\mathcal{H}}_{t}$.}
\item\label{lower-bound-reward}
From Assumption \ref{Reward+Multifunction}, it follows that for any $t\in \NN^{*}$,
$$
\ds \sup_{\mathbb{P}\in \boldsymbol{\mathcal{S}}} \int_{\boldsymbol{\mathcal{K}}_{t}} r_{t}^{+}(H_{t},A_{t}) d\mathbb{P}_{|\mathbf{H}_{t}\times \mathbf{A}_{t}} \leq \sup_{\mathbb{P}\in \boldsymbol{\mathcal{S}}} \int_{[K_{\epsilon} \times \mathbf{A}_{t}]\cap\boldsymbol{\mathcal{K}}_{t}} r_{t}^{+}(H_{t},A_{t})
d\mathbb{P}_{|\mathbf{H}_{t}\times \mathbf{A}_{t}} +\epsilon
$$
for some $\epsilon>0$ and a closed set $K_{\epsilon}$ in $\boldsymbol{\mathcal{H}}_{t}$. However, $r_{t}$ is bounded above on $[K_{\epsilon}\times  \mathbf{A}_{t}]\cap\boldsymbol{\mathcal{K}}_{t}$ and so,
$\ds \sup_{\mathbb{P}\in \boldsymbol{\mathcal{S}}} \int_{\boldsymbol{\mathcal{K}}_{t}} r_{t}^{+}(H_{t},A_{t}) d\mathbb{P}_{|\mathbf{H}_{t}\times \mathbf{A}_{t}} <\infty$.
\end{enumerate}
\end{remark}

\begin{proposition} \label{Reward-upper-semicompact}
Suppose Assumptions \ref{Regular-Hypotheses} and \ref{Reward+Multifunction} hold. The mapping $\mathcal{J}$ is upper semicompact.
\end{proposition}
\textbf{Proof:} We have to show that for any $\beta \in \RR$, the set $\{\mathbb{P}\in \boldsymbol{\mathcal{S}} : \mathcal{J}(P)\geq \beta\}$ is relatively compact in $\boldsymbol{\mathcal{S}}$ for the weak topology.
We will proceed along the line described in \cite[Lemma 4.2]{balder92} to show the result. The main difference is that in our case, the reward function is not necessarily bounded above and that the action sets are history-dependent leading to the introduction of a new set of conditions given by Assumption \ref{Reward+Multifunction}.

There exist compact sets $\widehat{\mathbf{X}}_{t}$ and $\widehat{\mathbf{A}}_{t}$ such that
$\mathbf{X}_{t}\in\mathfrak{B}(\widehat{\mathbf{X}}_{t})$ and $\mathbf{A}_{t}\in\mathfrak{B}(\widehat{\mathbf{A}}_{t})$ for $t\in \NN^{*}$. 
Equipped with the weak topology, $\boldsymbol{\mathcal{P}}(\mathbf{H}_{\infty})$ and $\boldsymbol{\mathcal{S}}$ are topological subspaces of the compact space $\boldsymbol{\mathcal{P}}(\widehat{\mathbf{H}}_{\infty})$ where 
$\widehat{\mathbf{H}}_{\infty}=\prod_{t\in \NN^{*}} \mathbf{X}_{t} \times \mathbf{A}_{t}$.
Consider $\{\mathbb{P}_{j}\}_{j\in \NN}$ a sequence in $\{\mathbb{P}\in \boldsymbol{\mathcal{S}} : \mathcal{J}(\mathbb{P})\geq \beta\} \subset \boldsymbol{\mathcal{P}}(\widehat{\mathbf{H}}_{\infty})$. Since 
$\boldsymbol{\mathcal{P}}(\widehat{\mathbf{H}}_{\infty})$ is compact, there exists a subsequence of $\{\mathbb{P}_{j}\}_{j\in \NN}$ (still denote by $\{\mathbb{P}_{j}\}_{j\in \NN}$) that converges to
$\mathbb{P}_{\infty}\in \boldsymbol{\mathcal{P}}(\widehat{\mathbf{H}}_{\infty})$. In order to get the the result, it is sufficient to show that $\mathbb{P}_{\infty}\in \boldsymbol{\mathcal{S}}$ or equivalently that
$\mathbb{P}_{\infty | \widehat{\mathbf{X}}_{1}}=\nu$ and for any $k\in \NN^{*}$
\begin{align}
\mathbb{P}_{\infty | \widehat{\mathbf{H}}_{k}\times \widehat{\mathbf{A}}_{k}}(\boldsymbol{\mathcal{K}}_{k})=1,
\label{Const-1} \\
\mathbb{P}_{\infty | \widehat{\mathbf{H}}_{k+1}}=P_{\infty | \widehat{\mathbf{H}}_{k}\times \widehat{\mathbf{A}}_{k}}Q_{k}.
\label{Const-2}
\end{align}
Clearly, we have $\mathbb{P}_{\infty | \widehat{\mathbf{X}}_{1}}=\nu$. The other two equalities will be shown by induction. Let us assume that (\ref{Const-1}) and
(\ref{Const-2}) hold for $k\in \{1,\ldots,t-1\}$.

Let us first show that $\mathbb{P}_{\infty | \widehat{\mathbf{H}}_{t} \times \widehat{\mathbf{A}}_{t}}(\boldsymbol{\mathcal{K}}_{t})=1$. To do so, it is sufficient to prove, for the function $v \colon \boldsymbol{\mathcal{H}}_{t}\times \widehat{\mathbf{A}}_{t} \to [-\infty,+\infty[$ defined by $v=r_{t}\wedge 0$ on $\boldsymbol{\mathcal{K}}_{t}$ and $v=-\infty$ otherwise, that
$$
\ds \int_{\boldsymbol{\mathcal{H}}_{t}\times \widehat{\mathbf{A}}_{t}} v d\mathbb{P}_{\infty | \widehat{\mathbf{H}}_{t}\times \widehat{\mathbf{A}}_{t}} > -\infty.
$$
Indeed, in such a case, the probability measure $\mathbb{P}_{\infty | \widehat{\mathbf{H}}_{t}\times \widehat{\mathbf{A}}_{t}}$ is necessarily supported by $\boldsymbol{\mathcal{K}}_{t}$.

For $\epsilon>0$, let us denote by $K_{\epsilon}$ the closed subset of $\boldsymbol{\mathcal{H}}_{t}$ satisfying Assumption \ref{Reward+Multifunction} and write $Z_{\epsilon}=[K_{\epsilon}\times \mathbf{A}_{t}]\cap\boldsymbol{\mathcal{K}}_{t}$.
By hypothesis, $r_{t | Z_{\epsilon}}$ is strongly coercive and so is $\big[r_{t | Z_{\epsilon}}\big]_{K_{\epsilon} \times \widehat{\mathbf{A}}_{t}}$ where
$$
\big[r_{t | Z_{\epsilon}}\big]_{K_{\epsilon} \times\widehat{\mathbf{A}}_{t}}=\left\{
\begin{array}{rl}
r_{t } &\text{on } Z_{\epsilon},\\
-\infty&\text{on } (K_{\epsilon}\times \widehat{\mathbf{A}}_{t})\setminus Z_{\epsilon}.
\end{array}
\right.
$$
Now, observe that $v_{| K_{\epsilon} \times  \widehat{\mathbf{A}}_{t}}=\big[r_{t | Z_{\epsilon}}\big]_{K_{\epsilon} \times \widehat{\mathbf{A}}_{t}} \wedge 0$. Thus, by items (i) and (iii) of Proposition 2.2 in \cite{balder92}, it follows that
$v_{| K_{\epsilon} \times  \widehat{\mathbf{A}}_{t}}$ is strongly coercive and so, upper semicontinuous.

By the induction hypothesis, we easily obtain  that
$\mathbb{P}_{\infty | \widehat{\mathbf{H}}_{t} \times \widehat{\mathbf{A}}_{t}}(\boldsymbol{\mathcal{H}}_{t} \times \widehat{\mathbf{A}}_{t})=1$.
Recalling that for any $j\in \NN^{*}$, $\mathbb{P}_{j | \widehat{\mathbf{H}}_{t}\times \widehat{\mathbf{A}}_{t}}(\boldsymbol{\mathcal{H}}_{t}\times \widehat{\mathbf{A}}_{t})=1$,
we get that the restriction of $\mathbb{P}_{j | \widehat{\mathbf{H}}_{t}\times \widehat{\mathbf{A}}_{t}}$ to $\boldsymbol{\mathcal{H}}_{t}\times \widehat{\mathbf{A}}_{t}$ converges weakly to
the restriction of $\mathbb{P}_{\infty | \widehat{\mathbf{H}}_{t}\times \widehat{\mathbf{A}}_{t}}$ to $\boldsymbol{\mathcal{H}}_{t}\times \widehat{\mathbf{A}}_{t}$.
Now, combining the Portmanteau theorem and \eqref{w-tight-condition}, it follows that $\ds \mathbb{P}_{\infty | \widehat{\mathbf{H}}_{t}\times \widehat{\mathbf{A}}_{t}}\big((\boldsymbol{\mathcal{H}}_{t}\setminus K_{\epsilon})\times \widehat{\mathbf{A}}_{t}\big)
\leq \liminf_{j\rightarrow\infty} \mathbb{P}_{j | \widehat{\mathbf{H}}_{t}\times \widehat{\mathbf{A}}_{t}}\big((\boldsymbol{\mathcal{H}}_{t}\setminus K_{\epsilon})\times \widehat{\mathbf{A}}_{t}\big) <\epsilon$
since $K_{\epsilon}$ is a closed subset of $\boldsymbol{\mathcal{H}}_{t}$.
Therefore, we can apply Corollary \ref{lem-tight-usc} to the function $v$ and the set of probability measures given by the restriction of
$\mathbb{P}_{j | \widehat{\mathbf{H}}_{t}\times \widehat{\mathbf{A}}_{t}}$ to $\boldsymbol{\mathcal{H}}_{t}\times \widehat{\mathbf{A}}_{t}$
with $j\in \NN\cup\{\infty\}$. We obtain $\ds \limsup_{j\rightarrow \infty} \int_{\boldsymbol{\mathcal{H}}_{t}\times \widehat{\mathbf{A}}_{t}} v d \mathbb{P}_{j | \widehat{\mathbf{H}}_{t}\times \widehat{\mathbf{A}}_{t}}
\leq \int_{\boldsymbol{\mathcal{H}}_{t}\times \widehat{\mathbf{A}}_{t}} v d\mathbb{P}_{\infty | \widehat{\mathbf{H}}_{t}\times \widehat{\mathbf{A}}_{t}}$.
However, recalling that $\mathcal{J}(\mathbb{P}_{j})\geq \beta$ we get with Assumption \ref{Condition-C} that $\beta\leq\mathcal{J}(\mathbb{P}_{j})\leq \mathcal{J}_{t}(\mathbb{P}_{j})+1$.
From Remark \ref{Hyp-discussion}.(\ref{lower-bound-reward}), it implies that $\ds \inf_{j\in \NN} \int_{\boldsymbol{\mathcal{H}}_{t}\times \widehat{\mathbf{A}}_{t}} (r_{t}\wedge 0) d\mathbb{P}_{j | \widehat{\mathbf{H}}_{t}\times \widehat{\mathbf{A}}_{t}} >-\infty$.
Therefore, $\ds \int_{\boldsymbol{\mathcal{H}}_{t}\times \widehat{\mathbf{A}}_{t}} v d\mathbb{P}_{\infty | \widehat{\mathbf{H}}_{t}\times\widehat{\mathbf{A}}_{t}} > -\infty$
showing that $\mathbb{P}_{\infty | \widehat{\mathbf{H}}_{t}\times \widehat{\mathbf{A}}_{t}}(\boldsymbol{\mathcal{K}}_{t})=1$ by definition of $v$, as required.

Now, following exactly the same arguments as in the second part of the proof of Theorem 2.1 in \cite{balder89} and using Assumption \ref{Continuity-transition-kernel},
we obtain that \eqref{Const-2} is valid for $k=t$, giving the last part of the result.
\hfill$\Box$

\begin{proposition} \label{Reward-upper-semicontinuous}
Suppose Assumptions \ref{Condition-C} and \ref{Reward+Multifunction} hold. The mapping $\mathcal{J}$ is upper semicontinuous.
\end{proposition}
\textbf{Proof:} Consider $t\in  \NN^{*}$ and $\epsilon>0$. According to Assumption \ref{Reward+Multifunction}  there exists $K_{\epsilon}$, a closed subset of $\boldsymbol{\mathcal{H}}_{t}$, such that
the restriction of $r_{t}$ to $Z_{\epsilon}=[K_{\epsilon}\times \mathbf{A}_{t}]\cap\boldsymbol{\mathcal{K}}_{t}$ is strongly coercive and bounded above and 
$\ds \sup_{\pi\in\Pi} \int_{\boldsymbol{\mathcal{K}}_{t}\setminus Z_{\epsilon}} \big[1\vee r_t^{+}(H_{t},A_{t}) \big] d\mathbb{P}^{\pi}_{| \mathbf{H}_{t} \times \mathbf{A}_{t}}<\epsilon$.
By item (i) of Proposition 2.2 in \cite{balder92}, the restriction of $r_{t}$ to $Z_{\epsilon}$ is upper semicontinuous.
Observe that $Z_{\epsilon}$ is a closed subset of $\boldsymbol{\mathcal{K}}_{t}$.
Therefore, from Corollary \ref{lem-tight-usc} we obtain that the mapping defined on $\boldsymbol{\mathcal{S}}$ by $\ds \mathbb{P}\to \int_{\boldsymbol{\mathcal{K}}_{t}} r_{t}(H_{t},A_{t}) d\mathbb{P}^{\pi}_{| \mathbf{H}_{t} \times \mathbf{A}_{t}}$ is upper semicontinuous and bounded above.
Now, taking into account Assumption \ref{Condition-C}, an application of Proposition 10.1 in \cite{schal75b} gives the result.
\hfill$\Box$\\

We are now able to state our main result.

\begin{theorem}\label{main-result}
Suppose that Assumptions \ref{Regular-Hypotheses} and \ref{Reward+Multifunction} hold and that there exists a strategic probability measure $P_m\in\mathbb{P}$ such that $\mathcal{J}(P_m)>-\infty$. Then there exists a policy $\pi^\ast\in \Pi$ such that
$$
\sup_{P\in \mathbb{P}}\mathcal{J}(P)=\mathcal{J}(P^{\pi^\ast}).
$$ 
\end{theorem}
\textbf{Proof:}
From Proposition \ref{Reward-upper-semicompact}, $\mathcal{J}$ is upper semicompact and so, the set $\{\mathcal{J}\geq \mathcal{J}(P_m)\}$ is relatively compact.
The map $\mathcal{J}$ being also upper semicontinuous according to Proposition \ref{Reward-upper-semicontinuous}, it admits a maximum on the compact set given by the closure of $\{\mathcal{J}\geq \mathcal{J}(P_m)\}$ and  the result follows.
\hfill$\Box$

\section{Examples}
\label{sec-examples}
This section provides examples illustrating our results.
The first example described a controlled model for which the reward function is unbounded and not strongly coercive on its domain of definition and takes the value $-\infty$. It is shown that this model satisfies our assumptions.

A set of hypotheses has been introduced in \cite{jaskiewicz11} to ensure in particular the existence of an optimal policy for Makov decision processes with unbounded rewards. It is shown in the second example that our conditions are satisfied {in such a setting. Moreover, we would like to emphasize that this set of conditions is satisfied for a large class of economical models as described in  \cite[Section 5]{jaskiewicz11}.
For the sake of completeness we describe one of such model at the end of the Section \ref{example-Nowak}.}

\subsection{An example with unbounded and non strongly coercive reward function}  
\label{example-general}
We consider a model with  state space $\mathbf{X}_{t}= [0,1]$ and action space $\mathbf{A}_{t}=\mathbb{N}^\ast$. 
Let us introduce $$\mathbf{A}(x)=\begin{cases} \{1,\ldots,p\} & \text{ if } x\in [0,1/2]\\
\mathbb{N}^\ast & \text{ if } x\in ]1/2,1], \end{cases}$$ for some $p\in\mathbb{N}^\ast$.
The sequence of multifunctions $\{\Psi_{t}\}_{t\in \NN^{*}}$ are defined recursively by setting $\Psi_{1} \colon  \mathbf{X}_{1} \to 2^{\mathbf{A_{1}}}  \backslash \{\emptyset\}$
with $\Psi_{1}(x)=\mathbf{A}(x)$ and  $\Psi_{t} : \boldsymbol{\mathcal{K}}_{t-1} \times \mathbf{X}_{t} \rightarrow 2^{\mathbf{A_{t}}} \backslash \{\emptyset\}$ for $t\geq 2$ 
given by $\Psi_{t}(x_{1},a_{1},\ldots,x_{t-1},a_{t-1},x_{t})=\mathbf{A}(x_{t})$ where $\boldsymbol{\mathcal{K}}_{t}\subset\mathbf{H}_{t}\times \mathbf{A}_{t}$ denotes the graph of $\Psi_{t}$.
The transition probability function from time $t$ to time $t+1$ is 
$$Q_t( dy | x_{1},a_{1},\ldots,x_{t-1},a_{t-1},x_{t},a_{t})=\beta(t+x_{t}/a_{t},5/2)(dy)$$ 
where $\beta(\alpha_{1},\alpha_{2})$ denotes the beta probability distribution on $[0,1]$ with parameters $(\alpha_{1},\alpha_{2})\in \RR_{+}^{2}$.
The reward functions we consider are given by
$$r_{t}(x_{1},a_{1},\ldots,x_{t-1},a_{t-1},x_{t},a_{t})= \begin{cases}
0 & \text{ if $(x_{t},a_{t}) \in \{0\}\times \NN^{*}$} \\
\frac{1}{a_{t}^{2}} \frac{1}{\sqrt{x_{t}}}I_{]0,1/2]}(x_{t})-a_{t} \frac{1}{\sqrt{1-x_{t}}}I_{]1/2,1[}(x_{t}) & \text{ if $(x_{t},a_{t}) \in ]0,1[\times \NN^{*}$} \\
-\infty, & \text{ if $(x_{t},a_{t}) \in \{1\} \times \NN^{*}$.}
\end{cases}$$
The initial probability measure $\nu$ is $\beta(1,2)$.

Observe that the reward function is unbounded, takes the value $-\infty$ and is not upper semicontinuous at point $(x_{1},a_{1},\ldots,x_{t-1},a_{t-1},0,a_{t})\in \boldsymbol{\mathcal{H}}_{t}$ and therefore not strongly coercive. Moreover,
the action sets available at each step to the decision maker are not compact.
{Let us show that Assumptions \ref{Regular-Hypotheses} and \ref{Reward+Multifunction} are satisfied. To check that Assumtpion \ref{Reward+Multifunction} holds, we will use the approach developed in Remark \ref{Hyp-discussion}(\ref{New-condition-1-sufficiency}). For any positive integer $t\geq1$ and $h_t=(x_{1},a_{1},\ldots,x_{t})\in \boldsymbol{\mathcal{H}}_{t}$, we set
\begin{equation}\label{eq:ex:maj}
\Phi_t(h_t)=I_{\{0\}}(x_t)+\frac{1}{\sqrt{x_t}}I_{]0,1/2]}(x_t)+I_{]1/2,1]}(x_t).
\end{equation}
This measurable function satisfies, for any $h_t\in \boldsymbol{\mathcal{H}}_{t}$,
\begin{eqnarray}
 \sup_{a_{t}\in \Psi_{t}(h_{t})} \big[1\vee r_{t}^{+}(h_{t},a_{t}) \big] \leq \Phi_{t}(h_{t}).
\end{eqnarray}
We have that for any positive integer $t\geq2$, $\epsilon\in [0,1/2]$, $\gamma\in [0,1/2]$ and any $(h_{t-1},a_{t-1})\in \boldsymbol{\mathcal{K}}_{t-1}$,
\begin{align}
&\int_{\mathbf{X}_{t}}\mathbf{I}_{\boldsymbol{\mathcal{K}}_{t-1}\times([0,\epsilon[\union]1/2,1/2+\gamma[)}(h_{t}) \Phi_{t}(h_{t}) Q_{t}(dx_{t}|h_{t-1},a_{t-1})\nonumber \\
& \leq  \int_{]0,\epsilon[} \frac{1}{\sqrt{y}} Q_{t}(dx_{t}|h_{t-1},a_{t-1})
+ Q_{t}(]1/2,1/2+\gamma [|h_{t-1},a_{t-1}) \nonumber\\
~&\leq \frac{1}{\mathsf{B}(t+x_{t-1}/a_{t-1},5/2)}
\bigg[  \int_{]0,\epsilon[} \frac{1}{\sqrt{y}} y^{t-1} dy  +  \int_{]1/2,1/2+\gamma [}  y^{t-1} dy \bigg] \nonumber\\
~&\leq \frac{1}{\mathsf{B}(t+1,5/2)} \Big[  \frac{1}{t-1/2} \epsilon^{t-1/2} + \frac{1}{t} \big[ (1/2+\gamma)^{t}-(1/2)^{t} \big]  \Big] \nonumber\\
& \leq  \frac{Kt^{7/2}}{t-1/2} \Big[ \epsilon^{t-1/2}+\big[ (1/2+\gamma)^{t}-(1/2)^{t} \big] \Big],
\label{Ineq-exemple2}
\end{align}
for some positive constant $K$ and where $\mathsf{B}(\alpha_{1},\alpha_{2})$ denotes the beta function with parameters $(\alpha_{1},\alpha_{2})\in \RR_{+}^{2}$.
For $t=1$, similar calculations lead to
\begin{equation}
\int_{\mathbf{X}_{1}}\mathbf{I}_{[0,\epsilon[\union]1/2,1/2+\gamma[}(x_1) \Phi_{1}(x_1)\nu(dx_1)\leq \nu(]1/2,1/2+\gamma [)+\frac{2}{\mathsf{B}(1,2)}\sqrt{\epsilon}\big(1-\epsilon/3\big).\label{Ineq-exemple2bis}
\end{equation}
}
Using equation \eqref{eq:ex:maj} and choosing $\epsilon=1/2$ and $\gamma=0$ in equation \eqref{Ineq-exemple2}, we can show that for $t\geq2$ and any policy $\pi\in \Pi$,
\begin{align*}
\mathbb{E}^\pi \big[ r^+_{t}(H_{t},A_t)\big] \leq \mathbb{E}^\pi \Big[ \mathbf{I}_{[0,\frac{1}{2}]}(X_t) \big[1\vee r^+_{t}(H_{t},A_t)\big]\Big]\leq \frac{Kt^{7/2}}{t-1/2} (1/2)^{t-1/2}\leq c t^{5/2} 2^{-t}
\end{align*}
for some positive constant $c$ and so, Assumption \ref{Condition-C} is satisfied since the series with main term $t^{5/2}2^{-t}$ is convergent.
Clearly, Assumption \ref{Continuity-transition-kernel} is satisfied.

Now, regarding Assumption \ref{Reward+Multifunction}, 
let us consider $t\in\NN^*$ and  $K_{\epsilon}$ defined by $\boldsymbol{\mathcal{K}}_{t-1} \times ([\epsilon,1/2]\cup[1/2+\epsilon,1])$
(with the slight abuse of notation $\boldsymbol{\mathcal{K}}_{0}\times ([\epsilon,1/2]\cup[1/2+\epsilon,1])$ means $[\epsilon,1/2]\cup[1/2+\epsilon,1]$) a closed subset of $\boldsymbol{\mathcal{H}}_{t}$ for any $\epsilon\in ]0,1/4]$.
By using equations \eqref{Ineq-exemple2} and \eqref{Ineq-exemple2bis}, the set $K_{\epsilon}$ satisfies, for $t\geq2$,
$$
\lim_{\epsilon\rightarrow 0} \sup_{(h_{t-1},a_{t-1})\in \boldsymbol{\mathcal{K}}_{t-1}}\int_{\mathbf{X}_{t}}\mathbf{I}_{\boldsymbol{\mathcal{H}}_{t}\setminus K_{\epsilon}}(h_{t}) \Phi_{t}(h_{t}) Q_{t}(dx_{t}|h_{t-1},a_{t-1})=0,
$$ 
and for $t=1$, using equation \eqref{Ineq-exemple2bis}
$$
\lim_{\epsilon\rightarrow 0} \int_{\mathbf{X}_{1}}\mathbf{I}_{[0,\epsilon[\union]1/2,1/2+\epsilon[}(x_1) \Phi_{1}(x_1)\nu(dx_1)=0.
$$
{This shows the existence of a closed subset of $\mathbf{X}_{1}$ (respectively, $\boldsymbol{\mathcal{H}}_{t+1}$)
satisfying condition \eqref{sufficient-tight-1-condition} (respectively, condition \eqref{sufficient-tight-2-condition}).
Now it remains to prove that the restriction of $r_{t}$ to $[K_{\epsilon} \times \mathbf{A}_{t}]\cap\boldsymbol{\mathcal{K}}_{t}$ is bounded above and strongly coercive to get the sufficient conditions proposed for Assumption 
 \ref{Reward+Multifunction}, in Remark \ref{Hyp-discussion}(\ref{New-condition-1-sufficiency}).}
Clearly, the restriction of $r_{t}$ to $[K_{\epsilon} \times \mathbf{A}_{t}]\cap\boldsymbol{\mathcal{K}}_{t}$ is bounded above. Let us show that it is strongly coercive.
Indeed, for $t\in\mathbb{N}^\ast$, let $\{g^{k},x^{k},a^{k}\}_{k\in \NN}$ be a sequence in $[\boldsymbol{\mathcal{K}}_{t-1} \times([\epsilon,1/2]\cup[1/2+\epsilon,1]) \times \mathbf{A}_{t}]\cap\boldsymbol{\mathcal{K}}_{t}$ such that $\{(g^{k},x^{k})\}_{k\in \NN}$ converges to $(g^\ast,x^{\ast})\in \boldsymbol{\mathcal{K}}_{t-1} \times ([\epsilon,1/2]\cup[1/2+\epsilon,1])$ as $k$ tends to infinity and $\ds \limsup_{k\rightarrow\infty} r_{t}(g^{k},x^{k},a^{k})>-\infty$.
Let us show that there necessarily exits a converging subsequence $\{a^{\phi(k)}\}_{k\in \NN}$ of $\{a^{k}\}_{k\in \NN}$ to $a^{*}$ such $a^{*} \in \mathbf{A}(x^{*})$ and 
$\ds \limsup_{k\rightarrow\infty} r_{t}(g^{k},x^{k},a^{k})\leq r_{t}(g^{*},x^{*},a^{*})$.

We first prove that there exists a subsequence $\{a^{\phi(k)}\}_{k\in \NN}$ of $\{a^k\}_{k\in \NN}$ such that $\{a^{\phi(k)}\}_{k\in \NN}$ converges towards some $a^\ast\in \mathbf{A}(x^\ast)=\mathbb{N}^\ast$.
Assume that $x^\ast$ is in $[\epsilon,1/2]$. Since $\{x^k\}$ is valued in $[\epsilon,1/2]\cup[1/2+\epsilon,1]$, this means that there is some $k_0$ such that for any $k\geq k_0$ we have $x_k\in[\epsilon,1/2]$ and thus
$$
\mathbf{A}(x^k)=\{1,\ldots,p\}.
$$
Therefore, $a_k$ is in the compact set $\{1,\ldots,p\}$ for all $k\geq k_0$ and has a convergent subsequences to some $a^\ast\in\{1,\ldots,p\}=\mathbf{A}(x^\ast)$.\\
Now, assume that $x^\ast$ is in $[1/2+\epsilon,1]$. Since $\{x^k\}_{k\in \NN}$ is valued in $[\epsilon,1/2]\cup[1/2+\epsilon,1]$, this means that there is some $k_0$ such that for any $k\geq k_0$ we have $x^k\in[1/2+\epsilon,1]$ and thus
$$
\mathbf{A}(x^k)=\mathbb{N}^\ast.
$$
We proceed now by contradiction. Assume that for any subsequences $\{a^{\phi(k)}\}_{k\in \NN}$ of $\{a^k\}_{k\in \NN}$, either $\{a^{\phi(k)}\}_{k\in \NN}$ diverges or it converges to some $a^\ast\notin \mathbf{A}(x^\ast)=\mathbb{N}^\ast$. Of course, the latter claim can not happen: if the integer valued $\{a^{\phi(k)}\}_{k\in \NN}$ converges to some $a^\ast$, then, necessarily, $a^\ast\in\mathbb{N}^\ast$.
Thus, assume that any subsequences $\{a^{\phi(k)}\}_{k\in \NN}$ of $\{a^k\}_{k\in \NN}$ diverges. This implies that the original sequence goes to infinity. Then, since for any $k\geq k_0$ we have $x_k\in[1/2+\epsilon,1]$,
$$
r_t(g^k,x^k,a^k)=-a^k\frac{1}{\sqrt{1-x^k}}.
$$
Therefore, $\ds \limsup_{k\rightarrow\infty} r(g^k,x^k,a^k)=-\infty$, showing a contradiction. Consequently, there exists a subsequence $\{a^{\phi(k)}\}_{k\in \NN}$ of $\{a^k\}_{k\in \NN}$ such that $\{a^{\phi(k)}\}_{k\in \NN}$ converges towards some $a^\ast\in \mathbf{A}(x^\ast)=\mathbb{N}^\ast$.

Now, in both cases, $x^\ast \in [\epsilon,1/2]$ or $x^\ast \in [1/2+\epsilon,1]$, the upper semi-continuity of $r_t$ on $[\boldsymbol{\mathcal{K}}_{t-1} \times([\epsilon,1/2]\cup[1/2+\epsilon,1]) \times \mathbf{A}_{t}$ implies that $\ds \limsup_{k\rightarrow\infty} r_{t}(g^{k},x^{k},a^{k})\leq r_{t}(g^{*},x^{*},a^{*})$.
Consequently, the restriction of $r_{t}$ to $[K_{\epsilon} \times \mathbf{A}_{t}]\cap\boldsymbol{\mathcal{K}}_{t}$ is strongly coercive showing that Assumption
\ref{Reward+Multifunction} is satisfied.

\subsection{A model by A. Ja\'{s}kiewicz and A. Nowak}
\label{example-Nowak}
The objective of this example is to show that our set of conditions are more general than those introduced in \cite{jaskiewicz11} to study discounted Markov decision processes with unbounded (from above and below)
reward function. In \cite{jaskiewicz11}, the authors consider the following discounted model with unbounded rewards where ${\cal X}$ (respectively, ${\cal A}$) is the Borel state (respectively, action) space.
The decision maker can choose the action in a compact set ${\cal A}(x)\subset {\cal A}$ when the state process is in $x\in {\cal X}$. The set valued mapping $\psi(x)={\cal A}(x)$ defined on ${\cal X}$ is assumed to be upper semicontinuous. The transition probability denoted by $q$ on ${\cal X}$ given the graph $\boldsymbol{\mathcal{K}}$ of $\psi$ is weakly continuous.
The reward function is given by $\beta^{t}u$ where the discount factor $\beta\in (0,1)$ and the reward function $u$ is upper semicontinuous on $\boldsymbol{\mathcal{K}}$. In addition to these classical hypotheses, it is assumed that there exists a sequence
$\{{\cal X}_{j}\}_{j\in \NN}$ of non-empty Borel subsets of ${\cal X}$ satisfying ${\cal X}=\cup_{j\in \NN} \mathring{{\cal X}}_{j}$ (where $\mathring{{\cal X}}_{j}$ denotes the interior of ${\cal X}_{j}$) and from the set ${\cal X}_{j}$ only transitions to ${\cal X}_{j+1}$ are allowed, in other words the transition kernel $q$ satisfies $q({\cal X}_{j+1} | x,a)=1$ for any $x\in {\cal X}_{j}$ and $a\in {\cal A}(x)$.
Moreover, the reward function $u$ is piecewise bounded on ${\cal X}$, that is $\sup_{x\in {\cal X}_{j}}\sup_{a\in {\cal A}(x)} u^{+}(x,a)=m_{j}<\infty$ with $\sum_{t\in \NN} \beta^{t}m_{j+t}<\infty$ for any $j\in \NN$.
By using the dynamic programming approach, it is proved in Theorem 1 in \cite{jaskiewicz11}, the existence of an optimal policy which is stationary maximizing $\mathbb{E}^{\pi}\big[\sum_{t\in \NN} \beta^{t}u(X_{t},A_{t})\big]$ over $\Pi$ when the initial distribution is degenerated to the point $x\in {\cal X}$.
It can be easily shown that the previous model can be embedded to our framework by choosing a set ${\cal X}_{j}$ satisfying $x\in {\cal X}_{j}$ and defining $\mathbf{X}_{t}={\cal X}_{j+t-1}$ for $t\in \NN^{*}$, $\mathbf{A}_{t}={\cal A}$ 
and $\Psi_{1} \colon  \mathbf{X}_{1} \to 2^{\mathbf{A}_{1}}  \backslash \{\emptyset\}$ with $\Psi_{1}(x)={\cal A}(x)$ for any $x\in \mathbf{X}_{1}$ and for $t\geq 2$,
$\Psi_{t} : \boldsymbol{\mathcal{K}}_{t-1} \times \mathbf{X}_{t} \rightarrow 2^{\mathbf{A}_{t}}$ with $\Psi_{t}(x_{1},a_{1},\ldots,x_{t-1},a_{t-1},x_{t})={\cal A}(x_{t})$ for any $(x_{1},a_{1},\ldots,x_{t-1},a_{t-1},x_{t})\in \boldsymbol{\mathcal{K}}_{t-1} \times \mathbf{X}_{t}$
and where $\boldsymbol{\mathcal{K}}_{t}\subset\mathbf{H}_{t} \times \mathbf{A}_{t}$ denotes the graph of $\Psi_{t}$. In our context, 
the stochastic kernel $Q_{t}$ on $\mathbf{X}_{t+1}$ given $\boldsymbol{\mathcal{K}}_{t}$ is defined by $Q_{t}(\cdot | x_{1},a_{1},\ldots,x_{t},a_{t}) =q(\cdot |x_{t},a_{t})$
and the reward function is given by $r_{t}:\boldsymbol{\mathcal{K}}_{t}\rightarrow [-\infty,+\infty[$ with
$r_{t}(x_{1},a_{1},\ldots,x_{t},a_{t})=\beta^{t}u(x_{t},a_{t})$ for any $(x_{1},a_{1},\ldots,x_{t},a_{t})\in \boldsymbol{\mathcal{K}}_{t}$.
Let us show that this model satisfies Assumptions \ref{Regular-Hypotheses} and \ref{Reward+Multifunction}.
Clearly, we have
$$\sup_{p\geq n} \sup_{\pi\in\Pi}  \bigg[ \sum_{t=n}^{p} \mathbb{E}^\pi \big[ r_t(X_{t},A_{t})\big] \bigg]^{+}\leq\sum_{t=n}^{\infty} \beta^{t}m_{j+t}$$
and so Assumption \ref{Condition-C} is satisfied since $\sum_{t\in \NN} \beta^{t}m_{j+t}<\infty$.
Moreover, Hypothesis \ref{Continuity-transition-kernel} holds since the stochastic kernel $q$ is weakly continuous.
Finally, it is easy to show that Assumption \ref{Reward+Multifunction} is satisfied. Indeed, $r_{t}$ is upper semicontinuous implying that $r_{t}$ is strongly continuous on the graph of $\Psi_{t}$ since $\Psi_{t}$ is upper semicontinuous with compact values. Now, by construction $r_{t}$ is bounded above on $\boldsymbol{\mathcal{K}}_{t}$ showing the claim. 

\bigskip

{To conclude this section, we present now an economic model borrowed from \cite[Section 5, Example 1]{jaskiewicz11}. We will show that Assumptions \ref{Regular-Hypotheses} and \ref{Reward+Multifunction} hold for this example.
Let ${\cal X}=[0,\infty)$ be the set of all possible capital stocks. The variable $X_t$ represents a capital stock at the beginning of the period $t$, during which a portion $A_t\in {\cal A}(X_t)=[0,X_t]$ of the capital is consumed. The evolution of the capital stock is given, for $t\in\mathbb{N}$, by
\begin{equation}\label{eq-ex2-ev}
X_{t+1}=(1+\rho)(X_t-A_t)+\Xi_t,
\end{equation}
with initial condition $X_1$. In the evolution equation \eqref{eq-ex2-ev}, $\rho>0$ is a constant rate of growth and $\Xi_t$ is a random income received in period $t$. The random variables $\{\Xi_t\}_{t\in\mathbb{N}^\ast}$ are assumed to be i.i.d. with probability distribution $\mu$ supported by $[0,z]$ for some $z\geq1$. The objective is to maximize the expected total utility of consumption, given by $r_t(X_t,A_t)=(A_t)^\sigma$ with $\sigma\in(0,1)$. 
For a positive real $d>0$ and $t\in\mathbb{N}^\ast$, introduce
$$
{\cal X}_t=[0,k_t],\quad {\cal A}={\cal X}.
$$
with $k_t=(1+\rho)^td+z(1+\rho)[(1+\rho)^{t-1}-1]/\rho$ and
$$
Q_t(C|x,a)=\int_0^z I_B((1+\rho)(x-a)+\xi)\mu(d\xi),
$$
where $C\in\mathfrak{B}({\cal X}_{t+1})$, $x\in {\cal X}_t$ and $a\in{\cal A}(x)$. Observe that the reward function $r_{t}$ is not bounded from above on ${\cal X}\times {\cal A}$. However,
it has been shown in \cite[Section 5, Example 1]{jaskiewicz11} that the set of hypotheses introduced in \cite{jaskiewicz11} are satisfied provided that
$m_t=k^\sigma_t$ and $\beta(1+\rho)^\sigma<1$. Therefore, by applying the arguments developed in the previous paragraph, we can claim that Assumptions \ref{Regular-Hypotheses} and \ref{Reward+Multifunction} are satisfied.
}

\appendix 

\section{Appendix : Balder's lemma}\label{app:balder}

We detail here the proof of \cite[Lemma 2.4]{balder89}, for the sake of clarity.
\begin{lemma}
\label{lem-balder}
Let $(\mu_{n})$ be a sequence of probability measures on a metric space $Y$ converging weakly to a probability measure $\mu$.
Consider a function $u:Y\rightarrow [-\infty,+\infty[$, bounded from above, satisfying the following conditions:
for any $\epsilon>0$,  there exists a closed subset $Y_{\epsilon}$ of $Y$ such that
\begin{align}
\label{Hyp-lem-balder}
\sup_{n\in\mathbb{N}}\mu_{n}(Y\setminus Y_{\epsilon}) <\epsilon
\end{align}
and the restriction of $u$ on $Y_{\epsilon}$ is upper semicontinuous.
Then, 
\begin{align}
\limsup_{n\rightarrow \infty} \int_{Y} u d\mu_{n} \leq \int_{Y} u d\mu.
\end{align}
\end{lemma}
\textbf{Proof:}
Let us assume in a first step that $u$ is also bounded from below. Let $\|u\|=\sup_{Y}|u|$ and define, for any $\epsilon>0$, the function
$$
u_\epsilon=u I_{Y_\epsilon} - \|u\|I_{Y\setminus Y_\epsilon}.
$$
Let us show that $u_\epsilon$ is upper semicontinuous on $Y$. For $\beta\in \mathbb{R}$, consider the level set
$$
A_\beta=\{ x\in Y : u_\epsilon(x) < \beta\}=\{x\in Y_\epsilon : u(x) < \beta\}\cup \{x\in Y\setminus Y_\epsilon : -\|u\| < \beta\}.
$$
Our aim is to show that $A_\beta$ is open. If $\beta\leq -\|u\|$, we clearly have $A_\beta=\emptyset$ which is an open set. Otherwise, we can write
$$
A_\beta=\{x\in Y_\epsilon : u_\epsilon(x)<\beta\}\cup (Y\setminus Y_\epsilon).
$$
Since $u$ is upper semicontinuous on $Y_\epsilon$, the level set $\{x\in Y_\epsilon : u_\epsilon(x) < \beta\}$ is open in $Y_\epsilon$, and so there exists an open set $O$ of $Y$ such that
$$
\{x\in Y_\epsilon : u_\epsilon(x)<\beta\}=Y_\epsilon\cap O.
$$
Thus $A_\beta=(Y_\epsilon\cap O)\cup (Y\setminus Y_\epsilon)$. Let $x\in A_\beta$. If $x\in Y\setminus Y_\epsilon$, by the closedness of $Y_\epsilon$, we can find $\eta>0$ such that $B(x,\eta)\subset  Y\setminus Y_\epsilon\subset A_\beta$. Otherwise, $x\in Y_\epsilon\cap O$. In this case, since $x\in O$ which is an open set, we can find $\eta'>0$ such that $B(x,\eta')\subset O$. Then
$$
B(x,\eta')\cap A_\beta= [B(x,\eta')\cap Y_\epsilon \cap O] \cup [ B(x,\eta')\cap (Y\setminus Y_\epsilon)]=[B(x,\eta')\cap Y_\epsilon] \cup [ B(x,\eta')\cap Y\setminus Y_\epsilon]=B(x,\eta').
$$ 
Thus $B(x,\eta')\subset A_\beta$ showing that $A_\beta$ is open. This implies that $u_\epsilon$ is upper semicontinuous on $Y$.\\
Remark that
$$
\sup_{n\in\mathbb{N}}\left|\int_Y u d\mu_n -\int_Y u_\epsilon d\mu_n\right|\leq 2\epsilon \|u\|.
$$
Now, using the fact that $u_\epsilon$ is upper semicontinuous and bounded on the whole space $Y$,
\begin{align*}
\limsup_{n\rightarrow \infty} \int_{Y} u d\mu_{n} \leq \limsup_{n\rightarrow \infty} \int_{Y} u_\epsilon d\mu_{n}+2\|u\|\epsilon 
\leq \int_{Y} u_\epsilon d\mu+2\|u\|\epsilon \leq \int_{Y} u d\mu+2\|u\|\epsilon,
\end{align*}
showing the result.\\
In the case where $u$ is no longer bounded from below, we introduce $u_m=u\vee (-m)$ for which the previous step holds. Then, we apply the monotone convergence theorem to obtain the result.
\hfill$\Box$


\end{document}